\documentstyle{amsart}

\newtheorem{Theorem}{Theorem}[section]
\newtheorem{Lemma}[Theorem]{Lemma}
\newtheorem{Corollary}[Theorem]{Corollary}

\newtheorem{Remark}[Theorem]{Remark}

\newtheorem{Example}[Theorem]{Example}

\begin{document}

\title{Generically finite morphisms}

\author{Steven Dale Cutkosky}
\thanks{  partially supported by NSF}

\maketitle
\section{introduction}

Suppose that $k$ is a field, and $f:Y\rightarrow X$ is a dominant, generically finite
morphism of complete $k$-varieties. If $Y$ and $X$ are complete curves, then
it is classical that $f$ is finite. If $Y$ and $X$ have dimensions $\ge 2$
$f$ need not be finite. The simplest example is  the blowup of a
nonsingular subvariety of a nonsingular projective variety.

It is however natural to ask the following question. 
Given a generically finite morphism
$f:Y\rightarrow X$ as above, does there exist a commutative diagram 
\begin{equation}\label{eqI1}
\begin{array}{lll}
Y_1&\stackrel{f_1}{\rightarrow}&X_1\\
\downarrow&& \downarrow\\
Y&\stackrel{f}{\rightarrow}&X
\end{array}
\end{equation}
such that $f_1$ is finite, $Y_1$ and $X_1$ are nonsingular complete $k$-varieties, and the vertical arrows are  birational?  The answer to this question is no, as is shown by
a theorem of Abhyankar (Theorem 11 \cite{Ab3}).  This theorem, (as shown in  Example \ref{Example2} of this paper) proves that such a diagram cannot always be constructed
even when $f:Y\rightarrow X$ is a $G$-equivariant morphism of complex projective surfaces,
where the extension of function fields $k(X)\rightarrow k(Y)$ is Galois with Galois
group $G$.

In the theory of resolution of singularities a modified version of this question
is important. 
\vskip .2truein
\noindent {\bf Question 1.} 
With $f:Y\rightarrow X$ as above, is it possible to construct a diagram (\ref{eqI1})
such that $f_1$ is finite, $Y_1$ and $X_1$ are  complete $k$-varieties such that $Y_1$ is nonsingular, $X_1$ is normal and the vertical arrows are  birational? 
\vskip .2truein
This question has
been posed by Abhyankar (with the further  conditions that $Y_1\rightarrow Y$ is a sequence
of blowups of nonsingular subvarieties and $Y$, $X$ are projective) explicitely on page 144 of \cite{Ab6}, where it is called the ``weak simultaneous resolution global conjecture'' and implicitely in the paper
\cite{Ab3}. 

As positive evidence for this conjecture, Abhyankar proves a local form of this conjecture
for 2 dimensional function fields over an algebraically closed field of arbitrary characteristic \cite{Ab1}, \cite{Ab4}, (this is the two dimensional case of the ``weak simultaneous resolution local  conjecture'' \cite{Ab6}). 

An important case where Question 1 has a positive answer 
is for  generically finite morphisms $f:Y\rightarrow X$ of projective varieties, over a field $k$ of characteristic
zero,  which induce a Galois
extension of function fields. We give a simple proof  in Theorem \ref{Theorem7}.
We can  construct (with this Galois assumption) a diagram (\ref{eqI1}) such that  the 
conclusions of Question 1 hold, and $X_1$ has normal toric singularities.
This is a relative version of Theorem 7, \cite{Ab3}.

In Theorem \ref{Example1} of this paper we give a counterexample to Question 1. 
The example is of a generically finite morphism $Y\rightarrow X$ of nonsingular, projective surfaces, defined over an algebraically closed field $k$ of characteristic not equal to 2.
This counterexample is necessarilly then a counterexample to the 
``weak simultaneous resolution global conjecture''.

As the ``weak simultaneous resolution local conjecture'', posed by Abhyankar on page 144 of \cite{Ab6} is true in characteristic 0 (we prove it in \cite{C3}
as a corollary to the local monomialization theorem, Theorem 1.1 \cite{C2}, and prove a stronger 
statement in Theorem \ref{Theorem1} of this paper),
Theorem \ref{Example1} also gives a counterexample to the philosophy that a theorem
which is true in valuation theory should also be true in the birational geometry of projective varieties. This is
the philosophy which led to successful proofs of resolution of singularities for surfaces
and 3-folds, in characteristic zero by Zariski (\cite{Z1}, \cite{Z5}), and in positive characteristic by Abhyankar (\cite{Ab1},\cite{Ab4},\cite{Ab5}). Recently there has been progress
on the important problem of local uniformization in positive characteristic in higher dimensions (c.f. \cite{HRW}, \cite{Ku},\cite{Moh}, \cite{T}). Ramification of morphisms of algebraic
surfaces in positive characteristic is analyzed in \cite{CP} and \cite{CP2}.

We  prove in Theorem \ref{Theorem8} that Question 1, and
the ``weak simultaneous resolution global conjecture''
are almost true, as it always possible (over fields of characteristic zero)  to construct a diagram (\ref{eqI1}) where 
$f_1:Y_1\rightarrow X_1$ is a quasi-finite morphism  of integral,   finite
type $k$-schemes, $Y_1$ is nonsingular, $X_1$ has normal toric singularities, the
vertical morphisms are birational and every $k$-valuation of $k(X)$ has a center on $X_1$, 
every $k$-valuation of $k(Y)$ has a center on $Y_1$. That is, the answer to Question 1 becomes true if we weaken the condition that the vertical arrows are proper by not insisting that these morphisms be separated.

The essential technical result used in the proof of 
Theorem \ref{Theorem8} is  the Local
Monomialization Theorem, Theorem 1.1 of \cite{C2}.  
Local monomialization is used to prove a strengthened version of the 
``weak simultaneous local conjecture'', Theorem \ref{Theorem1} of this paper, which allows us to
construct local solutions of the problem, which are patched in an arbitrary manner
(this is where separatedness is lost) to construct $X_1$ and $Y_1$. 

We will now give an overview of the proof of the construction of Theorem \ref{Example1},
which is the counterexample
to Question 1 and the ``weak simultaneous global conjecture''. We will use some of the notation
explained in the following section on notations.

Suppose that $K^*$ is a finite  extension of an algebraic function field $K$, defined
over an algebraically closed field $k$ of characteristic zero, $\nu^*$ is a $k$-valuation of
$K^*$ and $\nu$ is its restriction to $K$.  By Theorem \ref{Theorem1},  there exists an algebraic regular local ring $S$ with quotient field $K^*$ dominated by $\nu^*$ and an algebraic normal local
ring $R$ with quotient field $K$ such that $S$ lies over $R$ ($S$ is the localization at a maximal ideal  of the integral closure of $R$ in $K^*$). This can be refined \cite{CP2} to show that 
there exist $S$ and $R$ as above such that the quotient of value groups $\Gamma_{\nu^*}/\Gamma_{\nu}$
acts faithfully on the power series ring $\hat S$ by $k$-algebra automorphisms so that 
\begin{equation}\label{eqI2}
\hat R\otimes_{R/m_R}k'\cong (\hat S\otimes_{S/m_S}k')^{\Gamma_{\nu^*}/\Gamma_{\nu}}
\end{equation}
where $k'$ is an algebraic closure of $S/m_S$.
In some special cases, such as rational rank 2 valuations of algebraic function fields
of dimension two, where $[K:K^*]$ is not divisible by $\text{char}(k)$, this construction
is stable under quadratic transforms of $S$. We give a direct proof  in this paper.

Let $k$ be an algebraically closed field,  
$\overline \nu$ be the rational rank 2 valuation  on
$\overline K=k(u,v)$ and  $L_1$ be the  $q$-cyclic extension  of $\overline K$,
where $q$ is a prime distinct from $\text{char}(k)$, constructed in Theorem 11 \cite{Ab3} of Abhyankar (the construction is recalled in Theorem \ref{Theorem10} of this paper).
$\overline R=k[u,v]_{(u,v)}$ is dominated by $\overline \nu$.
The extension $\overline K\rightarrow L_1$ has the property that if $S_1$ is an algebraic regular local ring with 
quotient field $L_1$ which lies over
an algebraic normal  local ring $R_1$ with quotient field $\overline K$ such that $R_1$ is dominated
by  $\overline\nu$, and contains  $\overline R$,  then $R_1$ is singular. 

Examining this example, we see that there is a unique extension $\nu_1$ of
$\overline \nu$ to $L_1$, and
the quotient of value groups $\Gamma_{\nu_1}/\Gamma_{\overline\nu}\cong{\bold Z}_q$.
By (\ref{eqI2})
we have  $\hat R_1\cong\hat{S_1}^{{\bold Z}_q}$. Since $\hat R_1$ is singular
and $\hat S_1$ is a power series ring in two variables over $k$, the algebraic
fundamental group of $\hat R_1$ is
$$
\pi_1(\text{spec}(\hat R_1)-m_{\hat R_1})\cong{\bold Z}_q.
$$

We now consider the extension $\nu_2$ of $\overline\nu$ to a particular $p$-cyclic
extension $L_2$ of $\overline K$ where $p$ is a prime such that $p\ne q$ and $p\ne \text{char}(k)$.
We have $\Gamma_{\nu_2}/\Gamma_{\overline\nu}\cong {\bold Z}_p$, and if
$S_2$ is an algebraic regular local ring with quotient field $L_2$ which is dominated by $\nu_2$ which contains $\overline R$, and if there exists an algebraic normal local ring $R_2$
with quotient field $\overline K$ which lies below $S_2$,
then $\hat R_2\cong\hat{S_2}^{{\bold Z}_p}$ by (\ref{eqI2}). If $R_2$ is singular, the algebraic
fundamental group of $\hat R_2$ is then
$$
\pi_1(\text{spec}(\hat R_2)-m_{\hat R_2})\cong{\bold Z}_p.
$$

We then construct a morphism of projective nonsingular $k$-surfaces 
$\Phi:Y\rightarrow X$ such that $X$ has the function field $\overline K$ and $\overline R$ is the
local ring of a point on $X$. $Y$ is constructed in such a way that $Y$ splits into
two sheets over $\text{spec}(\overline R)$, and there are points on these two sheets which are formally the same as
extensions  $\overline R\rightarrow S_1$, $\overline R\rightarrow S_2$ into algebraic regular local rings with respective quotient fields $L_1$ and $L_2$ which are dominated by the respective valuations $\nu_1$ and $\nu_2$.  We then use these formal embeddings
to construct extensions $\overline\nu_1$ and $\overline \nu_2$ of $\nu$ to  the
function field $L_0$ of $Y$.

These extensions $\overline\nu_i$ have the property that if $Y_1$ is  nonsingular and $Y_1\rightarrow Y$ is proper birational (so that it can be factored by blowups of points)
then the map $Y_1\rightarrow X$ is formally isomorphic at the centers of the valuations
$\overline \nu_1$ and $\overline \nu_2$ to the extensions of $\overline R$ by the corresponding 
sequences of quadratic transforms of the local rings $S_1$ and $S_2$ along the 
respective valuations $\nu_1$ and $\nu_2$. 

Now suppose that we can construct  a diagram
$$
\begin{array}{lll}
Y_1&\rightarrow&X_1\\
\downarrow&&\downarrow\\
Y&\rightarrow&X
\end{array}
$$
such that $Y_1\rightarrow X_1$ is finite, $Y_1$ is nonsingular, $X_1$ is normal, and the
vertical arrows are proper and birational. Let $R_1$ be the local ring of the center of
$\overline \nu$ on $X_1$, $S(1)$ be the local ring of the center of $\overline\nu_1$ on $Y_1$, and let $S(2)$ be the local ring of the center of $\overline \nu_2$ on $Y_1$.
Since $\overline \nu_1$ and $\overline \nu_2$ both extend $\overline \nu$, 
and $\Gamma_{\overline\nu_1}/\Gamma_{\overline\nu}\cong{\bold Z_q}$,
$\Gamma_{\overline\nu_2}/\Gamma_{\overline\nu}\cong{\bold Z_p}$,
we must have that
$$
\widehat{S(1)}^{{\bold Z}_q}\cong\widehat R_1\cong \widehat{S(2)}^{{\bold Z}_p}.
$$
We then have that $R_1$ is singular, by our contruction of $\overline\nu_1$, and
thus the algebraic fundamental group $\pi_1(\text{spec}(\hat{R}_1)-m_{\hat R_1})$ has 
simultaneously order $p$ and order $q\ne p$
which is impossible.

\section{notations} 
 We will denote the
maximal ideal of a local ring $R$ by $m_R$.
We will denote the  quotient field of a domain $R$ by $QF(R)$.
Suppose that $R\subset S$ is an inclusion of local rings. We will say that $R$ dominates
$S$ if $m_S\cap R=m_R$.  
Suppose that $K$ is an algebraic function field over a field $k$. We will say that
a subring  $R$ of $K$ is algebraic if $R$ is essentially of finite type over $k$.
Suppose that $K^*$ is a finite extension of an algebraic function field $K$,
 $R$ is a local ring with $QF(K)$ and $S$ is a local ring with $QF(K^*)$. We will say
that $S$ lies over $R$ and $R$ lies below $S$ if $S$ is a localization at a maximal ideal of the integral closure of $R$ in $K^*$. 
If $R$ is a local ring, $\hat R$ will denote the completion of $R$ at its maximal ideal.

Good introductions to the valuation theory which we require in this paper can be found  in Chapter VI of \cite{ZS} and in \cite{Ab4}.  
A valuation $\nu$ of $K$ will be called a $k$-valuation if $\nu(k)=0$. 
We will denote by $V_{\nu}$ the associated valuation ring, which necessarily contains $k$. 
A valuation ring
$V$ of $K$ will be called a $k$-valuation ring if $k\subset V$.  The value group of a valuation $\nu$ will be denoted by $\Gamma_{\nu}$. If $X$ is an integral
$k$-scheme with function field $K$, then a point $p\in X$ is called a center of the
valuation $\nu$ (or the valuation ring $V_{\nu}$) if $V_{\nu}$ dominates ${\cal O}_{X,p}$.
If $R$ is a subring of $V_{\nu}$ then the center of $\nu$ (the center of $V_{\nu}$)
on $R$ is the prime ideal $R\cap m_{V_{\nu}}$. 

Suppose that $R$ is a local domain. A monoidal transform $R\rightarrow R_1$ is a 
birational extension of local domains such that $R_1=R[\frac{P}{x}]_m$ where $P$ is
a regular prime ideal of $R$, $0\ne x\in P$ and $m$ is a prime ideal of $R[\frac{P}{x}]$
such that $m\cap R=m_R$. $R\rightarrow R_1$ is called a quadratic transform  if $P=m_R$.

If $R$ is regular, and $R\rightarrow R_1$ is a monodial transform, then there exists a regular sustem of parameters $(x_1,\ldots, x_n)$ in
$R$ and $r\le n$ such that
$$
R_1=R\left[\frac{x_2}{x_1},\ldots,\frac{x_r}{x_1}\right]_m.
$$

Suppose that $\nu$ is a valuation of the quotient field $R$ with valuation ring $V_{\nu}$
which dominates $R$. Then $R\rightarrow R_1$ is a monoidal transform along $\nu$
(along $V_{\nu}$) if $\nu$ dominates $R_1$. 

We follow the notation of \cite{Ha}. In particular, we do not require that a scheme be separated.

\section{A counterexample to global weak simultaneous resolution}
In this section we construct the following example. This gives a counterexample
to Question 1 stated in the introduction, as well as to the
``weak simultaneous resolution global conjecture'' stated by Abhyankar explicitely on page 144
\cite{Ab6} and implicitely in the paper \cite{Ab3}. As the ``weak simultaneous resolution local conjecture'' is true in characteristic 0 (We prove it in  \cite{C3}, and prove a stronger
version in Theorem \ref{Theorem1} of this paper),
Theorem \ref{Example1} also gives a counterexample to the philosophy that a theorem
which is true in valuation theory should also be true in the birational geometry of projective varieties.

\begin{Theorem}\label{Example1}
Suppose that $k$ is an algebraically closed field of characteristic 0 or of odd
prime characteristic. Then there  exists  
 a generically finite morphism $\Phi:Y\rightarrow X$ of projective nonsingular 
$k$-surfaces such that there does not
exist a commutative diagram
$$
\begin{array}{lll}
Y_1&\stackrel{\Phi_1}{\rightarrow}&X_1\\
\downarrow&&\downarrow\\
Y&\stackrel{\Phi}{\rightarrow}&X
\end{array}
$$
where the vertical arrows are birational and proper, $Y_1$ is nonsingular, $X_1$ is normal, 
and $\Phi_1$ is finite.
\end{Theorem}

Throughout this section we will suppose that $k$ is an algebraically closed field.

\begin{Lemma}\label{Lemma3}
Suppose that $L$ is a 2 dimensional algebraic function field over $k$. Suppose that $R$
is an algebraic regular local ring with quotient field $L$ and maximal ideal $m_{R}=(u,v)$.
Suppose that $\overline\nu$ is a rank 1, rational rank 2 valuation of $L$ such that
$\overline \nu$ dominates $R$,
$\overline\nu(u),\overline\nu(v)>0$ and $\overline \nu(u),\overline \nu(v)$ are
rationally independent. Then
\begin{enumerate}
\item The value group of $\overline\nu$ is $\Gamma_{\overline \nu}={\bold Z}\overline\nu(u)+{\bold Z}\overline\nu(v)$.
\item Suppose that $R\rightarrow R_1$ is a sequence of quadratic
transforms along $\overline \nu$. Then there exist regular parameters $(u_1,v_1)$ in 
$R_1$ and $a,b,c,d\in{\bold N}$ such that
$$
\begin{array}{ll}
u&=u_1^av_1^b\\
v&=u_1^cv_1^d
\end{array}
$$
with $ad-bc=\pm 1$.
\item There exists a unique extension $\hat\nu$ of $\overline \nu$ to $\hat L=QF(\hat R)$ which dominates $\hat{R}$. The value group of $\hat \nu$ is $\Gamma_{\hat\nu}=\Gamma_{\nu}$.
\item If $\nu_1$ is a valuation such that $\nu_1$ is equivalent to $\overline \nu$,
(and the value groups $\Gamma_{\overline\nu}$ and $\Gamma_{\nu_1}$ are embedded as subgroups
of ${\bold R}$) then 
$$
\frac{\nu_1(v)}{\nu_1(u)}=\frac{\overline \nu(v)}{\overline\nu(u)}.
$$
\end{enumerate}
\end{Lemma}

\begin{pf} 
Proof of 1. $f\in R$ implies there is an expression $f=\sum_{i+j=r}^{n-1}a_{ij}u^iv^j+h$ with
$a_{ij}\in k$, 
$r=\text{ord}(f)$, $h\in (m_{R})^n$, where $n$ is such that 
$$
n\overline\nu(m_{R})>\nu(\sum_{i+j=r}a_{ij}u^iv^j).
$$
Thus since $\overline \nu(u)$ and $\overline \nu(v)$ are rationally independent,
$$
\overline \nu(f)=\overline\nu(\sum_{i+j=r}^{n-1}a_{ij}u^iv^j)=\text{min}\{\overline \nu(u^iv^j)\mid r\le i+j\le n-1,\,\,\,a_{ij}\ne 0\}.
$$
Thus $\Gamma_{\overline\nu}={\bold Z}\overline\nu(u)+{\bold Z}\overline\nu(v)$.

Proof of 2. It suffices to prove this for a single quadratic transform.
We either have that $\overline \nu(u)>\overline \nu(v)$ or $\overline \nu(v)>\overline \nu(u)$. In the first case have that 
$$
R_1=R[\frac{u}{v},v]_{(\frac{u}{v},v)}
$$
and $\overline \nu(\frac{u}{v})$, $\overline\nu(v)$ are linearly independent over ${\bold Q}$.
In the second case we have that
$$
R_1=R[u,\frac{v}{u}]_{(u,\frac{v}{u})}
$$
and $\overline\nu(u)$, $\overline\nu(\frac{v}{u})$ are linearly independent over $\bold Q$.

Proof of 3. Define an extension $\hat\nu$ of $\overline\nu$ to $\hat L$ by
$$
\hat\nu(f)=\text{min}\{i\overline\nu(u)+j\overline\nu(v)\mid a_{ij}\ne0\}
$$
if $f\in\hat{R}$, and $f$ has the expression $f=\sum a_{ij}x^iy^j$ with $a_{ij}\in k$.  $\hat\nu$ is a valuation since
for $i,j,\alpha,\beta\in {\bold N}$, 
$$
i\overline\nu(u)+j\overline\nu(v)=\alpha\overline\nu(u)+\beta\overline\nu(v)
$$
implies $i=\alpha, j=\beta$.
$\hat \nu$ dominates
$\hat{R}$ and $\Gamma_{\hat\nu}=\Gamma_{\nu}$.

Suppose that $\tilde\nu$ is an extension of $\overline\nu$ to $\hat L$ which dominates
$\hat R$. 
Suppose that $f\in\hat R$. Write
$$
f=\sum_{i+j=r}^{\infty}a_{ij}u^iv^j,
$$
where $r=\text{ord}(f)$, $a_{ij}\in k$. There exists $n$ such that 
$$
\overline \nu(\sum_{i+j=r}a_{ij}u^iv^j)<n\overline\nu(m_R).
$$
Write
$$
f=\sum_{i+j=r}^{n-1}a_{ij}u^iv^j+g,
$$
with $g\in m_R^n\hat R$.
$$
\tilde \nu(\sum_{i+j=r}^{n-1}a_{ij}u^iv^j)=\overline\nu(\sum_{i+j=r}^{n-1}a_{ij}u^iv^j)
=\text{min}\{\overline\nu(u^iv^j)\mid a_{ij}\ne 0, r\le i+j\le n-1\}
<n\overline\nu(m_R)
$$
and $\tilde \nu(g)\ge n\overline\nu(m_R)$ so that
$$
\tilde\nu(f)=\text{min}\{\overline\nu(u^iv^j)\mid a_{ij}\ne 0\}=\hat\nu(f).
$$

Proof of 4.  As on page 653 \cite{Z1}, we consider the convergent factions $\frac{f_p}{g_p}$
of $\tau=\frac{\overline\nu(v)}{\overline\nu(u)}$. Set 
$$
\epsilon=f_{p-1}g_p-f_pg_{p-1}=\pm1.
$$
$\epsilon,-\tau+\frac{f_{p-1}}{g_{p-1}},\tau-\frac{f_p}{g_p}$ have the same signs.

For arbitrary $p$, we can define $u_1,v_1\in L$ by 
$$
u=u_1^{g_p}v_1^{g_{p-1}},
v=u_1^{f_p}v_1^{f_{p-1}}.
$$
$$
u_1^{\epsilon}=\frac{u^{f_{p-1}}}{v^{g_{p-1}}},
v_1^{\epsilon}=\frac{v^{g_p}}{u^{f_p}}.
$$
$$
\begin{array}{ll}
\epsilon \overline \nu(u_1)&=f_{p-1}\overline\nu(u)-g_{p-1}\overline\nu(v)\\
&=g_{p-1}\overline\nu(u)\left[\frac{f_{p-1}}{g_{p-1}}-\tau\right]
\end{array}
$$
which implies that $\overline\nu(u_1)>0$.
$$
\begin{array}{ll}
\epsilon\overline\nu(v_1)&=g_p\overline\nu(v)-f_p\overline\nu(u)\\
&=\overline\nu(u)g_p\left[ \tau-\frac{f_p}{g_p}\right]
\end{array}
$$
which implies $\overline \nu(v_1)>0$. Thus $\nu_1(u_1), \nu_1(v_1)>0$.
$$
f_{p-1}\nu_1(u)-g_{p-1}\nu_1(v)=\epsilon\nu_1(u_1),\,\,\,
-f_p\nu_1(u)+g_p\nu_1(v)=\epsilon\nu_1(v_1)
$$
imply
$$
\frac{f_{p-1}}{g_{p-1}}>\frac{\nu_1(v)}{\nu_1(u)}>\frac{f_p}{g_p}
$$
if $\epsilon=1$,
$$
\frac{f_{p-1}}{g_{p-1}}<\frac{\nu_1(v)}{\nu_1(u)}<\frac{f_p}{g_p}
$$
if $\epsilon=-1$. Since this holds for all $p$, 
$$
\frac{\nu_1(v)}{\nu_1(u)}=\frac{\overline \nu(v)}{\overline\nu(u)}.
$$
\end{pf}

\begin{Lemma}\label{Lemma7}
Suppose that $L$ is a 2 dimensional algebraic function field over $k$. Suppose that $R$
is an algebraic regular local ring with quotient field $L$ and maximal ideal $m_{R}=(u,v)$.
Suppose that $\nu_1$ is a rank 1, rational rank 2 valuation of $\hat L=QF(\hat R)$
such that $\nu_1(u)$, $\nu_1(v)>0$ are rationally independent and which dominates $\hat{R}$. Then $\overline \nu=\nu_1\mid L$ is a rank 1, rational
rank 2 valuation such that
$$
\Gamma_{\nu_1}={\bold Z}\overline \nu(u)+{\bold Z}\overline \nu(v)=\Gamma_{\overline \nu}.
$$
\end{Lemma}

\begin{pf} By arguments as in the proof of 3. of Lemma \ref{Lemma3}, we see that if $f=\sum_{i+j=r}^{\infty}a_{ij}u^iv^j\in\hat R$ with $a_{ij}\in k$,
then $\nu_1(f)=\text{min}\{\nu_1(u^iv^j)\mid a_{ij}\ne 0\}$. Thus 
$$
\Gamma_{\nu_1} ={\bold Z}\overline\nu(u)+{\bold Z}\overline\nu(v)=\Gamma_{\overline\nu}.
$$
\end{pf}

\begin{Remark} Suppose that $R$ is an algebraic regular local ring with quotient field $K$.
There exist many extensions of a given valuation $\nu$ of $K$ which
dominates $R$ to $QF(\hat R)$ which do not dominate $\hat R$. Let $K=k(x)$,
$\nu$ be the rank 1 discrete valuation  with valuation ring $V_{\nu}=k[x]_{(x)}$ such that $\nu(x)=1$. Set
$R=k[x]_{(x)}$.

Choose $f_1,\ldots,f_n\in \hat R=k[[x]]$ such that $x,f_1,\ldots,f_n$ are algebraically
independent over $k$, and choose $\gamma_1,\ldots,\gamma_n\in {\bold R}$ such that
$1,\gamma_1,\ldots,\gamma_n$ are linearly independent over ${\bold Q}$.
$K_1=k(x,f_1,\ldots,f_n)$ is a rational function field in $n+1$ variables, so we can extend $\nu$ to a rank 1, rational rank $n+1$ valuation $\nu_1$ of $K_1$ by setting $\nu_1(f_i)=\gamma_i$, $1\le i\le n$. By Proposition 2.22 \cite{Ab4} or Theorem 5', Section 4, Chapter VI \cite{ZS}, $\nu_1$
extends (up to equivalence) to a valuation $\hat \nu$ of $QF(\hat R)$ which we can normalize so that it is an extension of $\nu$.

Write $f_i=x^{m_i}\lambda_i$ where $\lambda_i\in\hat R$ is a unit series.
$\nu_1(\lambda_i)=\gamma_i-m_i\ne 0$. Since $\lambda_i$ and $\lambda_i^{-1}\in\hat R$,
$\hat R$ contains elements of negative $\hat\nu$ value, and thus $\hat\nu$ does not
dominate $\hat R$.
\end{Remark}

\begin{Lemma}\label{Lemma4}Suppose that $K\rightarrow K^*$ is a finite
extension of algebraic function fields over $k$ of dimension 2. Suppose that $\nu$ is a rank 1
rational rank 2 valuation of $K$, $\nu^*$ is an extension of $\nu$ to $K^*$. Suppose that
$R_0$ is an algebraic regular local ring with quotient field $K$, maximal ideal $m_{R_0}=(u,v)$, 
$S$ is an algebraic regular local ring with quotient field $K^*$,  maximal ideal $m_S=(x,y)$ and such that $S$ dominates $R$,
$$
\begin{array}{ll}
u&=x^ay^b\delta_1\\
v&=x^cy^d\delta_2
\end{array}
$$
for some natural numbers $a,b,c,d$ and units $\delta_1,\delta_2\in S$, and 
such that the characteristic of $k$ does not divide $ad-bc$.

Suppose that $V_{\nu^*}$ dominates $S$ and $\nu(u),\nu(v)$ are rationally independent over $\bold
Q$. Then
$$
\Gamma_{\nu}={\bold Z}\nu(u)+{\bold Z}\nu(v)
$$
and $\nu^*$ is a rank 1, rational rank 2 valuation of $K^*$ such that $\nu^*(x),\nu^*(y)$ are
rationally independent over $\bold Q$, and
$$
\Gamma_{\nu^*}={\bold Z}\nu^*(x)+{\bold Z}\nu^*(y).
$$
Suppose that $S\rightarrow S_1$ is a sequence of quadratic transforms along $\nu^*$.
Then $S_1$ has regular parameters $(\overline x_1,\overline y_1)$ such that
$$
\begin{array}{ll}
x&=\overline x_1^{\overline a}\overline y_1^{\overline b}\\ 
y&=\overline x_1^{\overline c}\overline y_1^{\overline d}
\end{array}
$$
with $\overline a\overline d-\overline b\overline c=\pm1$, and there exists a
(unique) algebraic regular
local ring $R_1$ with quotient field $K$ which lies below $S_1$. 
$\hat R_1\cong \hat{S_1}^{\Gamma_{\nu^*}/\Gamma_{\nu}}$, where
$\Gamma_{\nu^*}/\Gamma_{\nu}$ acts faithfully on $\hat{S_1}$ by $k$-algebra automorphisms, by multiplication
of $\overline x_1,\overline y_1$ by roots of unity in $k$.
\end{Lemma}

\begin{pf} $\nu^*$ has rational rank 2 and rank 1 since $K^*$ is finite over $K$
(Lemmas 1 and 2 of Section 11, Chapter VI \cite{ZS}). $\nu^*(x),\nu^*(y)$ are linearly
independent over $\bold Q$, so Lemma \ref{Lemma3} applies to $\nu$ and to $\nu^*$. 
We have an expression in $S_1$
$$
\begin{array}{ll}
u&=\overline x_1^{\tilde a}\overline y_1^{\tilde b}\tilde \delta_1\\
v&=\overline x_1^{\tilde c}\overline y_1^{\tilde d}\tilde \delta_2
\end{array}
$$
with natural numbers $\tilde a,\tilde b,\tilde c,\tilde d$ and units
$\tilde\delta_1$ and $\tilde\delta_2$ in $S_1$ such that the characteristic of $k$ does not divide $\tilde a\tilde d-\tilde b\tilde c$. 
Let 
$$
A=\left(\begin{array}{ll} \tilde a&\tilde b\\
\tilde c&\tilde d\end{array}\right),
$$
$d=\mid\tilde a\tilde d-\tilde b\tilde c\mid$. 
There exist regular parameters $\tilde x_1,\tilde y_1$ in $\hat S_1$ such that 
$$
u=\tilde x_1^{\tilde a}\tilde y_1^{\tilde b},\,\,\,\,\,
v=\tilde x_1^{\tilde c}\tilde y_1^{\tilde d}.
$$
Let $\omega$ be a $d^{th}$ root of unity in $k$. ${\bold Z}^2/A{\bold Z}^2$ acts
faithfully on $\hat S_1$ by $k$-algebra automorphisms. To $c\in {\bold Z}^2/A{\bold Z}^2$
the corresponding  $k$-algebra automorphism $\sigma_c$ of $\hat S_1$ is defined by
$$
\sigma_c(\tilde x_1) =\omega^{<B_1,c>}\tilde x_1\,\,\,\,\,
\sigma_c(\tilde y_1)=\omega^{<B_2,c>}\tilde y_1
$$
where $B_i$ is the $i^{th}$ row of $dA^{-1}=\pm\text{adj}(A)$. Since $k(u,v)\rightarrow
k(\tilde x_1,\tilde y_1)$ is Galois with Galois group ${\bold Z}^2/A{\bold Z}^2$, it follows
that $\hat S_1^{{\bold Z}^2/A{\bold Z}^2}$ is the completion of a $k$-algebra generated by
rational monomials $u^{\alpha_1}v^{\beta_1},\cdots,u^{\alpha_r}v^{\beta_r}$ (with
$\alpha_i,\beta_i\in{\bold Z}$ for all $i$).
$$
R_0[u^{\alpha_1}v^{\beta_1},\ldots,u^{\alpha_r}v^{\beta_r}]\subset \hat S_1\cap K=S_1.
$$
Let $R_1=R_0[u^{\alpha_1}v^{\beta_1},\ldots,u^{\alpha_r}v^{\beta_r}]_p$ where
$p=R_0[u^{\alpha_1}v^{\beta_1},\ldots,u^{\alpha_r}v^{\beta_r}]\cap m_{S_1}$.
$\hat R_1=\hat S_1^{{\bold Z}^2/A{\bold Z}^2}$ is normal, so
$R_1=\hat R_1\cap K$ is normal. Since $\sqrt{m_{R_1}S_1}=m_{S_1}$, $R_1$ lies below $S_1$
by Zariski's main Theorem (10.9 \cite{Ab5}). Uniqueness follows since  the condition
$R_1$ lies below $S_1$ implies $R_1=S_1\cap K$ by Proposition 1 (iv) \cite{Ab1}.
\end{pf}

\begin{Remark}The conclusion $\hat R_1\cong \hat{S_1}^{\Gamma_{\nu^*}/\Gamma_{\nu}}$
in Lemma \ref{Lemma4} is a special case of
a general result on ramification of valuations \cite{CP2}.
\end{Remark}

\begin{Lemma}\label{Lemma5} Suppose that $p$ is a prime such that $p$ is not
the characteristic of $k$ and ${\bold Z}_p$
acts diagonally and faithfully on the powerseries ring $k[[x,y]]$. Set
$R=k[[x,y]]^{{\bold Z}_p}$. Then $R$ is a normal local ring such that either
\begin{enumerate}
\item $R$ is regular and the algebraic fundamental group 
$$
\pi_1(\text{spec}(R)-m_R)=0
$$
or
\item $R$ is not regular, $R\rightarrow k[[x,y]]$ is unramified away from $m_R$ and the algebraic fundamental group
$$
\pi_1(\text{spec}(R)-m_R)\cong {\bold Z}_p.
$$
\end{enumerate}
\end{Lemma}

\begin{pf} Let $\omega$ be a primitive $p^{th}$ root of unity in $k$, $\sigma$ a generator of
${\bold Z}_p$. There exist integers $a,b$ with $0\le a,b<p$ such that
$$
\sigma(x)=\omega^ax,\,\,\,\,\,\sigma(y)=\omega^by.
$$

Suppose that $a=0$. Then $R=k[[x,y^p]]$ is regular. If $b=0$ then $R=k[[x^p,y]]$ is regular.
In both cases, 
$$
\pi_1(\text{spec}(R)-m_R)=\pi_1(\text{spec}(R))=\pi_1(k)=0
$$
 by the purity of the branch locus
(Theorems X 3.4, X 1.1  \cite{SGA2}).

Suppose that $a,b\ne 0$. Then for $1\le i\le p-1$ there exists a unique $j_i$ such 
that $bj_i\equiv ai\text{ mod }p$ with $0<j_i<p$. This implies that $x^{p-i}y^{j_i}$
is an invariant. Note that there exists an invariant of the form $x^{p-i_1}y$ for some
$0<i_1<p$, so that $j_{i_1}=1$. We will show that
$$
R=k[[x^p,x^{p-1}y^{j_1},\ldots,xy^{j_{p-1}},y^p]].
$$
We must show that any invariant monomial in $x$ and $y$ is a product of powers of these
$p+1$ monomials.

Suppose that $x^iy^j$ is invariant. Then  $ai+bj\equiv 0\text{ mod }p$. Write
 $i=\overline i+\lambda p$, $j=\overline j+\tau p$, with 
$0\le \overline i<p$, $0\le \overline j<p$.
$$
x^iy^j=x^{\overline i}y^{\overline j}x^{\lambda p}y^{\tau p}
$$
 $b\overline j\equiv -a\overline i\text{ mod }p$
implies $\overline i=\overline j=0$ or $\overline j=j_{p-\overline i}$.

Consider the finite map of normal local rings
$$
\Phi:Y=\text{spec}(k[[x,y]])\rightarrow X=\text{spec}(R).
$$
The ramification locus of $\Phi$ is defined by the $2\times 2$ minors of
$$
J(\Phi)=
\left(  \begin{array}{ll}
\frac{\partial(x^p)}{\partial x}& \frac{\partial(x^p)}{\partial y}\\
\frac{\partial (x^{p-1}y^{j_1})}{\partial x}&\frac{\partial(x^{p-1}y^{j_1})}{\partial y}\\
\vdots&\vdots\\
\frac{\partial(y^p)}{\partial x}&\frac{\partial(y^p)}{\partial y}
\end{array}\right).
$$
$$
\text{Det}\left(\begin{array}{ll}
y^{j_{p-1}}&j_{p-1}xy^{j_{p-1}-1}\\
 0&py^{p-1}\\
\end{array}\right)=py^{p-1+j_{p-1}}
$$
and
$$
\text{Det}\left(\begin{array}{ll} px^{p-1}&0\\
(p-i_1)x^{p-i_1-1}y&x^{p-i_1}\end{array}\right)=px^{2p-1-i_1}
$$
implies  $\sqrt{I_2(J(\Phi))}=(x,y)$. Thus $\Phi$ is 
unramified (and \'etale) away from $m_R$. 

Suppose that $S$ is a complete normal local domain such that $S$ is finite over $R$, and $R\rightarrow S$ is \'etale away from $m_R$. Let $T$ be the normalization of the image of
$S\otimes_Rk[[x,y]]$ in $QF(S)\otimes_{QF(R)}QF(k[[x,y,]])$. $k[[x,y]]\rightarrow T$
is \'etale away from $(x,y)$, so by the purity of the branch locus, and since $k$ is
algebraically closed, $\text{spec}(T)$ is a disjoint union of copies of $\text{spec}(k[x,y]])$.
A choice of one of these copies gives a factorization
$$
\text{spec}(k[[x,y]])\rightarrow\text{spec}(S)\rightarrow\text{spec}(R).
$$
Thus $\pi_1(X-m_R)\cong {\bold Z}_p$.
\end{pf}

Abhyankar constructs an example which  shows that we cannot in general take  $R$ to be regular
in general in Corollary \ref{Corollary1} (and thus  we cannot take $R$ to be regular in Theorem \ref{Theorem1}).

\begin{Theorem}\label{Theorem10}(Abhyankar)
There exists a two dimensional algebraic regular local ring $\overline R$ with quotient
field $\overline K$, a valuation $\overline \nu$ of $\overline K$ which dominates $\overline R$, and a finite extension $L_1$ of $\overline K$ such that if $\overline R_1$ is an algebraic regular local ring with quotient field $\overline K$ such that $\overline R\subset \overline R_1$ and $V_{\overline \nu}$ dominates $\overline R_1$, then there
is a unique normal algebraic local ring $\overline S$ with quotient field $L_1$ lying over $\overline R_1$.  $\overline S$ is not regular.
\end{Theorem}

\begin{pf} We give an outline of the construction, refering to Theorem 11 \cite{Ab3}
for details.

Let $\overline K=k(u,v)$ be a rational function field in two variables. Let $q>3$ be a prime such that
$q\ne \text{char}(k)$. Set $a=q-4$, 
Set
$$
\tau=a+\frac{1}{1+\frac{1}{a+\frac{1}{1+\frac{1}{a+\cdots}}}}\in{\bold R}-{\bold Q}.
$$
Define a rank 1, rational rank 2 valuation $\overline \nu$ on $\overline K$ by setting
$\overline \nu(u)=\tau$, $\overline\nu(v)=1$.

Set $\overline R=k[u,v]_{(u,v)}\subset V_{\overline\nu}$. 
Let
$$
L_1=\overline K[\overline z]/\overline z^q-uv^2.
$$
Let $z$ be the image of $\overline z$ in $L_1$.

Abhyankar shows that if $\overline R_1$ is an algebraic regular local ring with quotient
field $\overline K$ such that $\overline R\subset \overline R_1$, and $V_{\overline\nu}$ dominates $\overline R_1$, then there exists a unique normal algebraic local ring $\overline S$ with quotient field $L_1$ lying over $\overline R_1$
and $\overline S$ is not regular.
\end{pf}

By Lemma \ref{Lemma3}, $\Gamma_{\overline\nu}={\bold Z}+\tau{\bold Z}$. We will show that
there is a unique extension $\nu_1$ of $\overline \nu$ to $L_1$. First suppose that $\nu_1$ is 
a valuation of $L_1$ such that $V_{\nu_1}\cap \overline K=V_{\overline\nu}$.
Since $\nu_1$ must have rank 1 and rational rank 2 (Lemmas 1 and 2, Section 11, Chapter VI \cite{ZS}), we can assume that the value group of
$\nu_1$ is a subgroup of ${\bold R}$. We can then assume that
$\nu_1$ is normalized so that $\nu_1(v)=1$.  Since $\nu_1\mid \overline K$ is equivalent to $\overline\nu$, and $\nu_1(v)=1$, we have $\nu_1\mid \overline K=\overline\nu$ by Lemma \ref{Lemma3}. Thus $\nu_1$ is an
extension of $\overline\nu$. Since $\nu_1(z)=\frac{1}{q}(2+\tau)$, we have that
$\Gamma_{\nu_1}/\Gamma_{\overline\nu}\cong{\bold Z}_q$, and $\nu_1$ is the unique extension
of $\overline\nu$ to $L_1$, by corollary to Theorem 25, Section 12, Chapter VI \cite{ZS} and
Lemma 2.18 \cite{Ab4}.

Let $p$  be another prime such that $p\ne q$ and $p\ne\text{char}(k)$, and set
$$
L_2=\overline K[\overline w]/\overline w^p-uv^2.
$$
Let $w$ be the image of $\overline w$ in $L_2$.
By the same analysis as for $\nu_1$,
there is a unique extension $\nu_2$ of $\overline \nu$ to $L_2$. $\nu_2(w)=
\frac{1}{p}(2+\tau)$ and $\Gamma_{\nu_2}/\Gamma_{\overline\nu}\cong {\bold Z}_p$.

We remark that 
\begin{equation}\label{eq13}
\tau=q-4+\epsilon\text{ with }0<\epsilon<1.
\end{equation}

Set $x_1(1)=\frac{v}{z}$, $y_1(1)=\frac{z^2}{v}\in L_1$.
$$
\begin{array}{ll}
\nu_1(x_1(1))&=\nu_1(v)-\nu_1(z)=\\
&=\frac{q-2-\tau}{q}=\frac{2-\epsilon}{q}>0
\end{array}
$$
$$
\begin{array}{ll}
\nu_1(y_1(1))&=2\nu_1(z)-\nu_1(v)=\frac{2}{q}(2+\tau)-1\\
&=\frac{q-4+2\epsilon}{q}>0.
\end{array}
$$

Set $S_1=k[x_1(1),y_1(1)]_{(x_1(1),y_1(1))}$. $QF(S_1)=L_1$.
$\overline R\subset S_1\subset V_{\nu_1}$. 
\begin{equation}\label{eq18}
\begin{array}{ll}
u&=x_1(1)^{q-4}y_1(1)^{q-2}\\
v&=x_1(1)^2y_1(1).
\end{array}
\end{equation}
$\Gamma_{\nu_1}=\nu_1(x_1(1)){\bold Z}+\nu_1(y_1(1)){\bold Z}$.

We will now impose the further condition that $5\le q<p<2q-4$.
For example, we could take $q=11, p=13$ or $q=17, p=23$.  Set
$$
x_1(2)=\frac{v}{w},\,\,\, y_1(2)=\frac{w^2}{v}\in L_2.
$$ 
$$
\begin{array}{ll}
\nu_2(x_1(2))&=\nu_2(v)-\nu_2(w)\\
&=\frac{p-2-\tau}{p}=\frac{(p-q)+2-\epsilon}{p}>0
\end{array}
$$
$$
\begin{array}{ll}
\nu_2(y_1(2))&=2\nu_2(w)-\nu_2(v)=\frac{2}{p}(2+\tau)-1\\
&=\frac{2q-p-4+2\epsilon}{p}>0.
\end{array}
$$
Set 
$$
S_2=k[x_1(2),y_1(2)]_{(x_1(2),y_1(2)}.
$$
$QF(S_2)=L_2$. $\overline R\subset S_2\subset V_{\nu_2}$. 
\begin{equation}\label{eq19}
\begin{array}{ll}
u&=x_1(2)^{p-4}y_1(2)^{p-2}\\
v&=x_1(2)^2y_1(2).
\end{array}
\end{equation}

$\Gamma_{\nu_2}=\nu_2(x_1(2)){\bold Z}+\nu_2(y_1(2)){\bold Z}$.

We will now assume that $\text{char}(k)\ne 2$.

Let  $k[x,y,z_1]$ be a polynomial ring in $x,y,z_1$,
$$
f=z_1^2-1+x^my^n\in k[x,y,z_1]
$$
with $m,n$ odd and sufficiently large, as will be determined below. We will also assume that $m,n$ are not divisible by $\text{char}(k)$.
Then $f$ is irreducible. Set
$S_0=k[x,y,z_1]/(f)$. By abuse of notation, we will from now on identify $x,y,z_1$ with their
equivalence classes  in $S_0$. $S_0$ is smooth over $k$.
Suppose that $a_1,b_1,c_1,d_1\in{\bold N}$ are such that $a_1d_1-b_1c_1$ is not divisible by $\text{char}(k)$
and $a_2,b_2,c_2,d_2\in{\bold N}$ are such that $a_2d_2-b_2c_2$ is not divisible by $\text{char}(k)$. We now impose the conditions
$$
m>\text{max}\{\mid a_1-a_2\mid,\mid c_1-c_2\mid\}
$$
and
$$
n>\text{max}\{\mid b_1-b_2\mid,\mid d_1-d_2\mid\}.
$$
Let $R=k[u,v]$ be a polynomial ring in two variables.

 We define a $k$-algebra homomorphism
$$
R\rightarrow S_0
$$
by 
$$
\begin{array}{l}
u=x^{a_1}y^{b_1}(1-z_1)+x^{a_2}y^{b_2}(1+z_1)\\
v=x^{c_1}y^{d_1}(1-z_1)+x^{c_2}y^{d_2}(1+z_1)
\end{array}
$$
Consider the prime ideals $P_1=(x,y,z_1+1)$ and $P_2=(x,y,z_1-1)$ in $S_0$. 
In the local ring $(S_0)_{P_1}$ we have
$(P_1)_{P_1}=(x,y)$
since 
$$
(z_1+1)=-(z_1-1)^{-1}x^my^n.
$$ 
In the local ring $(S_0)_{P_2}$ we have $(P_2)_{P_2}=(x,y)$.

In $(S_0)_{P_1}$ we have
$$
\begin{array}{ll}
u&=x^{a_1}y^{b_1}(1-z_1-(z_1-1)^{-1}x^{a_2+m-a_1}y^{b_2+n-b_1})\\
&=x^{a_1}y^{b_1}\delta_1\\
v&=x^{c_1}y^{d_1}(1-z_1-(z_1-1)^{-1}x^{c_2+m-c_1}y^{d_2+n-d_1})\\
&=x^{c_1}y^{d_1}\delta_2
\end{array}
$$
where $\delta_1,\delta_2\in (S_0)_{P_1}$ are units.

Since $a_1d_1-b_1c_1$ is not divisible by $\text{char}(k)$, we have regular parameters $x_1(1),y_1(1)\in\widehat{(S_0)_{P_1}}$ such that 
\begin{equation}\label{eq16}
\begin{array}{ll}
u&=x_1(1)^{a_1}y_1(1)^{b_1}\\
v&=x_1(1)^{c_1}y_1(1)^{d_1}
\end{array}.
\end{equation}
This implies that $R\subset S_0$ is an inclusion.

By a similar calculation, we have  units $\epsilon_1,\epsilon_2\in (S_0)_{P_2}$
such that
$$
\begin{array}{ll}
u&=x^{a_2}y^{b_2}\epsilon_1\\
v&=x^{c_2}y^{d_2}\epsilon_2
\end{array}
$$
and 
regular parameters $x_1(2),y_1(2)\in\widehat{(S_0)_{P_2}}$ such that 
\begin{equation}\label{eq17}
\begin{array}{ll}
u&=x_1(2)^{a_2}y_1(2)^{b_2}\\
v&=x_1(2)^{c_2}y_1(2)^{d_2}
\end{array}.
\end{equation}

With the notation introduced in Theorem \ref{Theorem10}, we have $\overline R=R_{(u,v)}$
and $\overline K=QF(R)$.
Set $L_0=QF(S_0)$. $L_0$ is finite over $\overline K$  since $\overline K\rightarrow L_0$ is an inclusion of algebraic function fields of dimension 2.

In this construction, set 
$$
a_1=q-4,b_1=q-2, c_1=2, d_1=1
$$
 and 
$$
a_2=p-4,b_2=p-2,c_2=2,d_2=1,
$$
where $p,q$ are  the primes chosen 
in Theorem \ref{Theorem10}, and in the paragraph following Theorem \ref{Theorem10}.

Let $\Phi:Y\rightarrow X$ be a  morphism of smooth projective surfaces over $k$
which extends our map
$\text{spec}(S_0)\rightarrow \text{spec}(R)$. Such a map exists by resolution of singularities
of surfaces in characteristic $\ge 0$ (\cite{Ab5}, \cite{H2},  \cite{L3}). For $i=1,2$ we have commutative diagrams: 
\begin{equation}\label{eq14}
\begin{array}{lll}
\overline R=R_{(u,v)}&\rightarrow&(S_0)_{P_i}\\
\downarrow&&\downarrow\\
\hat{\overline R}=k[[u,v]]&\rightarrow &\widehat{(S_0)_{P_i}}=k[[x_1(i),y_1(i)]]
\end{array}
\end{equation}
with (by (\ref{eq16}) and (\ref{eq17}))
$$
\begin{array}{ll}
u&=x_1(i)^{a_i}y_1(i)^{b_i}\\
v&=x_1(i)^{c_i}y_1(i)^{d_i}.
\end{array}
$$
We further have commutative diagrams: 
\begin{equation}\label{eq15}
\begin{array}{lll}
\overline R&\rightarrow&S_i\\
\downarrow&&\downarrow\\
\hat{\overline R}=k[[u,v]]&\rightarrow&\hat S_i=k[[x_1(i),y_1(i)]]
\end{array}
\end{equation}
with (by  (\ref{eq18}) and (\ref{eq19}))
$$
\begin{array}{ll}
u&=x_1(i)^{a_i}y_1(i)^{b_i}\\
v&=x_1(i)^{c_i}y_1(i)^{d_i}.
\end{array}
$$
Diagrams (\ref{eq14}) and (\ref{eq15}) patch to give commutative diagrams for $i=1,2$ 
\begin{equation}\label{eq20}
\begin{array}{lll}
R_{(u,v)}&\rightarrow &(S_0)_{P_i}\\
\downarrow&&\downarrow\\
\hat{\overline R}&\rightarrow &\hat S_i
\end{array}
\end{equation}
Lemma \ref{Lemma3} implies that for $i=1,2$, there exists a unique extension $\hat \nu_i$ of $\nu_i$ to 
$QF(\hat S_i)$ which dominates $\hat S_i$ and $\Gamma_{\hat\nu_i}=\Gamma_{\nu_i}$.

For $i=1,2$, let $\overline\nu_i=\hat\nu_i\mid L_0$ (under the inclusion $L_0\subset QF(\hat S_i)$
induced by (\ref{eq20})).
$\Gamma_{\overline\nu_i}\cong \Gamma_{\nu_i}$ by Lemma \ref{Lemma7}.
$\overline \nu_i\mid \overline K=\overline\nu$ for $i=1,2$ where $\overline\nu$ is the
valuation of $\overline K$ introduced in Theorem \ref{Theorem10}.

Suppose that there exists a diagram
$$
\begin{array}{lll}
Y_1&\stackrel{\Phi_1}{\rightarrow}&X_1\\
\downarrow&&\downarrow\\
Y&\stackrel{\Phi}{\rightarrow}&X
\end{array}
$$
such that the vertical arrows are birational and proper, $Y_1$ is nonsingular, $X_1$ is normal, 
and $\Phi_1$ is finite. Then $Y_1\rightarrow Y$ is a sequence of blowups of points (Theorem II.1.1 \cite{Z6}). There exist commutative diagrams 

\begin{equation}\label{eq21}
\begin{array}{lll}
R_1&\rightarrow&S(1)\\
\uparrow&&\uparrow \lambda_1\\
\overline R&\rightarrow&(S_0)_{P_1}
\end{array}
\end{equation}
and 
\begin{equation}\label{eq22}
\begin{array}{lll}
R_1&\rightarrow&S(2)\\
\uparrow&&\uparrow \lambda_2\\
\overline R&\rightarrow&(S_0)_{P_2}
\end{array}
\end{equation}

where $R_1$ is the center of $\overline\nu$ on $X_1$, $S(1)$ is the local ring of the center of $\overline\nu_1$ on $Y_1$ and $S(2)$ is the local ring of the center of $\overline\nu_2$ on $Y_1$. $\lambda_1$ and
$\lambda_2$ are products of quadratic transforms.

By Lemma \ref{Lemma4} 
\begin{equation}\label{eq23}
\hat R_1\cong \widehat{S(1)}^{\Gamma_{\overline\nu_1}/\Gamma_{\overline \nu}}.
\end{equation}
and 
\begin{equation}\label{eq24}
\hat R_1\cong \widehat{S(2)}^{\Gamma_{\overline\nu_2}/\Gamma_{\overline\nu}}.
\end{equation}

By Theorem \ref{Theorem10}, (\ref{eq20}) with $i=1$ and (\ref{eq21}), $R_1$ is not regular. 
$$
\Gamma_{\overline\nu_1}/\Gamma_{\overline \nu}\cong {\bold Z}_q\text{ and }
\Gamma_{\overline\nu_2}/\Gamma_{\overline \nu}\cong {\bold Z}_p
$$
by our construction.

By Lemma \ref{Lemma5} and (\ref{eq23}), 
$$
\pi_1(\text{spec}(\hat R_1)-m_{\hat R_1})
\cong{\bold Z}_q
$$
 and by (\ref{eq24})
$$
\pi_1(\text{spec}(\hat R_1)-m_{\hat R_1})
\cong{\bold Z}_p.
$$
But $p\ne q$, so we have a contradiction.

\section{Ramification of valuations in algebraic function fields}

\begin{Theorem}\label{TheoremA} (Monomialization; Theorem 1.1 \cite{C2})
 Let $k$ be a field of characteristic zero, $K$ an algebraic function field over $k$, $K^*$ a finite algebraic extension of $K$,
$\nu^*$ a $k$-valuation of $K^*$. Suppose that $S^*$ is an algebraic  regular local ring with
quotient field $K^*$ which is dominated by $\nu^*$ and $R^*$ is an algebraic regular local ring with
quotient field $K$ which is dominated by $S^*$.
  Then there
exist sequences of  monoidal transforms 
 $R^* \rightarrow R_0$ and $S^* \rightarrow S$
such that $\nu^*$ dominates $S$, $S$ dominates $R_0$ and 
there are regular parameters $(x_1, .... ,x_n)$
in $R_0$,  $(y_1, ... ,y_n)$ in $S$, units $\delta_1,\ldots,\delta_n\in S$ and a matrix $A=(a_{ij})$ of
nonnegative integers such that  $\text{det}(A) \ne 0$ and 
 \begin{equation}\label{eq1*}
\begin{array}{lll}
x_1 &=& y_1^{a_{11}} \cdots y_n^{a_{1n}}\delta_1\\
&&\vdots\\   
x_n &=& y_1^{a_{n1}} \cdots y_n^{a_{nn}}\delta_n.
\end{array}
\end{equation}
\end{Theorem}

The standard theorems on resolution of singularities allow one to easily find $R_0$ and $S$ such that (\ref{eq1*}) holds,
but, in general, we will not  have the essential condition $\text{det}(a_{ij}) \ne 0$. The difficulty in the proof of this Theorem is to achieve the condition $\text{det}(a_{ij})\ne0$.

Let $\alpha_i$ be the images of $\delta_i$ in $S/m_{S}$ for $1\le i\le n$. Let 
$C=(a_{ij})^{-1}$, a matrix with rational coefficients. Define regular parameters 
$(\overline y_1,\ldots,\overline y_n)$ in $\hat{S}$ by
$$
\overline y_i=\left(\frac{\delta_1}{\alpha_1}\right)^{c_{i1}}\cdots\left(\frac{\delta_n}{\alpha_n}\right)^{c_{in}}y_i
$$
for $1\le i\le n$. We thus have relations 
\begin{equation}\label{eq25}
x_i=\alpha_i\overline y_1^{a_{i1}}\cdots\overline y_n^{a_{in}}
\end{equation}
with $\alpha_i\in S/m_S$ for $1\le i\le n$ in
$$
\hat R_0=R_0/m_{R_0}[[x_1,\ldots,x_n]]\rightarrow\hat S=S/m_S[[\overline y_1,\ldots,\overline y_n]].
$$

\begin{Remark}\label{RemarkT} Suppose that $k'$ is a field,
$A=(a_{ij})$ is an $n\times n$ matrix of natural numbers with $\text{det}(A)\ne 0$,
and we have an inclusion of lattices
$$
N'={\bold Z}^n\stackrel{A}{\rightarrow} N={\bold Z}^n
$$
with a corresponding inclusion of dual lattices
$$
M=\text{Hom}(N,{\bold Z})\rightarrow M'=\text{Hom}(N',{\bold Z}).
$$
Let $\sigma$ be the cone generated by the rows of $A$ in $N\otimes {\bold R}$,
$\sigma'$ be the cone generated by ${\bold N}^n$ in $N'\otimes{\bold R}$.
Suppose that $\hat\sigma\cap M$ is generated by $e_1,\ldots,e_r$. By the theory of
toric varieties (c.f. page 34 \cite{Fu}) we have inclusions of $k'$-algebras
$$
k'[\overline x_1,\ldots,\overline x_n]\rightarrow k'[\hat\sigma\cap M]=k'[\overline x^{e_1},
\ldots,\overline x^{e_r}]\rightarrow k'[\hat\sigma'\cap M']=k'[\overline y_1,\ldots,\overline y_n]
$$
with $\overline x_i=\overline y_1^{a_{i1}}\cdots\overline y_n^{a_{in}}$ for $1\le i\le n$.
$k'[\hat\sigma\cap M]$ is a normal ring with quotient field $k'(\overline x_1,\ldots,\overline x_n)$ and $k'[\hat\sigma'\cap M']$ is finite over $k[\hat\sigma\cap M]$.
\end{Remark}

If $k'={\bold C}$, this can be expressed in a particularly nice way.
$N/N'\cong {\bold Z}^n/A{\bold Z}^n$ acts on ${\bold C}[\overline y_1,\ldots,\overline y_n]$
by associating to $c\in N/N'$ the ${\bold C}$-algebra automorphism $\sigma_c$ defined by
$$
\sigma_c(\overline y_i)=\text{exp}^{2\pi i<F_i,c>}\overline y_i
$$
for $1\le i\le n$, where 
$$
A^{-1}=\left(\begin{array}{l} F_1\\ \vdots\\ F_n\end{array} \right).
$$

\begin{Theorem}\label{Theorem1} Let $k$ be a field of characteristic zero, $K$ an algebraic function field over $k$, $K^*$ a finite algebraic extension of $K$,
$\nu^*$ a $k$-valuation of $K^*$. Suppose that $S^*$ is an algebraic  local ring with
quotient field $K^*$ which is dominated by $\nu^*$ and $R^*$ is an algebraic local ring with
quotient field $K$ which is dominated by $S^*$.  Then there exists a commutative diagram  
\begin{equation}\label{eq3}
\begin{array}{lllllll}
R_0&\rightarrow& R&\rightarrow& S&\subset&V_{\nu^*}\\
\uparrow&&&&\uparrow\\
R^*&&\rightarrow&&S^*
\end{array}
\end{equation}
where $S^*\rightarrow S$  and $R^*\rightarrow R_0$ are sequences of monodial transforms along $\nu^*$ such that $R_0\rightarrow S$ have regular parameters of the form of the conclusions of 
Theorem \ref{TheoremA}, $R$ is an algebraic normal local ring with toric singularities, which is the localization of the blowup of an ideal in $R_0$, and
the regular local ring $S$ is the localization at a maximal
ideal of the integral closure of $R$ in $K^*$. 
\end{Theorem}

\begin{pf} 
By resolution of singularities \cite{H1} (c.f. Theorem 2.6, Theorem 2.9 \cite{C2}), we first reduce to the case where $R^*$ and $S^*$ are regular, and then
construct, by  the local monomialization theorem, Theorem \ref{TheoremA} a sequence of monodial transforms along $\nu^*$ 
\begin{equation}\label{eq6}
\begin{array}{lllll}
R_0&\rightarrow &S&\subset &V_{\nu^*}\\
\uparrow&&\uparrow\\
R^*&\rightarrow &S^*
\end{array}
\end{equation}
so that $R_0$ is a regular local ring with  regular parameters $(x_1,\ldots, x_n)$,
$S$ is a regular local ring with regular parameters $(y_1,\ldots, y_n)$, there are units
$\delta_1,\ldots,\delta_n$ in $S$, and a matrix of natural numbers
$A=(a_{ij})$ with nonzero determinant $d$ such that
$$
x_i=\delta_iy_1^{a_{i1}}\cdots y_n^{a_{in}}
$$
for $1\le i\le n$.

Let $k'$ be an algebraic closure of $S/m_S$. With the notation of (\ref{eq25}), set
$\overline x_i=\frac{x_i}{\alpha_i}$,
so that
$$
R_0\otimes_{R_0/m_{R_0}}k'\cong k'[[\overline x_1,\ldots,\overline x_n]]
\rightarrow S\otimes_{S/m_S}k'\cong k'[[\overline y_1,\ldots,\overline y_n]]
$$
is defined by 
$$
\overline x_i=\overline y_1^{a_{i1}}\cdots\overline y_n^{a_{in}},
$$
$1\le i\le n$.  With the notation of Remark \ref{RemarkT}, set $R=R_0[x^{e_1},\ldots,x^{e_r}]_m$
where $m=(x^{e_1},\ldots,x^{e_r})$. 
$$
R_0[x^{e_1},\ldots,x^{e_r}]\subset \hat S\cap K^*=S
$$
(by Lemma 2 \cite{Ab1}) and $m\subset m_{\hat S}\cap K^*=m_S$, so $S$ dominates $R$.
$\hat R=R_0/m_{R_0}[[\overline x^{e_1},\ldots,x^{e_r}]]$ is integrally closed in its
quotient field (by Remark \ref{RemarkT} and Theorem 32, Section 13, Chapter VIII \cite{Z2}),
so $R=\hat R\cap K$ is integrally closed in $K$. 
After possibly reindexing the $y_i$, we may assume that $d=\text{det}(A)>0$. Let $(b_{ij})$
be the adjoint matrix of $A$. Then
$$
\prod_{j=1}^nx_j^{b_{ij}}=\left(\prod_{j=1}^n\delta_j^{b_{ij}}\right)y_i^d\in R.
$$
Thus $\sqrt{m_RS}=m_S$,
so $R$ lies below $S$ by Zariski's Main Theorem (10.9 \cite{Ab5}).

\end{pf}

As an immediate consequence, we obtain a proof in characteristic zero of the ``weak
simultaneous resolution local conjecture''. which is stated explicitely on page 144 of
\cite{Ab6}, and is implicit in \cite{Ab3}. Abhyankar proves this for algebraic function
fields of dimension two and any characteristic in \cite{Ab1} and \cite{Ab4}.
In the paper \cite{C3}, we have given a direct proof of this result, also
as a consequence of Theorem \ref{TheoremA} (Theorem 1.1 \cite{C2}).

\begin{Corollary}\label{Corollary1}(Corollary1)(Theorem 1.1 \cite{C3}) Let $k$ be a field of characteristic zero, $K$ an algebraic function field over $k$, $K^*$ a finite algebraic extension of $K$,
$\nu^*$ a $k$-valuation of $K^*$, and $S^*$ an algebraic regular local ring with
quotient field $K^*$ which is dominated by $\nu^*$. Then for some sequence of monodial transforms $S^*\rightarrow S$ along $\nu^*$, there exists a normal algebraic local ring $R$ with quotient field $K$, such that the regular local ring $S$ is the localization at a maximal
ideal of the integral closure of $R$ in $K^*$.
\end{Corollary}

\begin{pf}  There exists a normal algebraic local ring $R^*$ with quotient field $K$ such that $\nu^*$ dominates $R^*$ (take $R^*$ to be the local ring of the center of $\nu^*$ on a
normal projective model of $K$).
There exists a finite type $k$-algebra $T$ such that the integral closure of $R^*$ in $K^*$ is a localization of $T$, and $T$ is generated over $k$ by $g_1,\ldots,g_m\in V^*=V_{\nu^*}$ such that
$\nu^*(g_i)\ge 0$ for all $i$. There exists a sequence of monoidal transforms
$S^*\rightarrow S_1$ along $\nu^*$ such that $T\subset S_1$ (Theorem 2.7 \cite{C2}).
$S_1$ dominates $R^*$.
After  replacing $S^*$ with $S_1$, we can assume that $S^*$ dominates $R^*$.
Theorem \ref{Theorem1} applies to this situation, so we can construct a diagram of the form (\ref{eq3}).
\end{pf}

When $K^*$ is Galois over $K$,
it is not difficult to construct using Galois theory and resolution of singularities
a regular local ring $S$ with quotient field $K^*$ and a normal local ring $R$ with
quotient field $K$ such that $S$ lies over $R$ (Theorem 7 \cite{Ab3}, Theorem \ref{Theorem7}), although even in the Galois case the full statements of Theorem  \ref{Theorem1} and Corollary \ref{Corollary1}
do not follow from these results (Theorem 7 \cite{Ab3}, Theorem \ref{Theorem7}). The general case of non Galois extensions is much more subtle, and not as well behaved, as can be seen from Theorem \ref{Example1}.

\section{Generically finite morphisms}

Suppose that $f:Y\rightarrow X$ is a dominant, generically finite morphism of complete
$k$-varieties. In this section we construct a commutative diagram
$$
\begin{array}{lll}
\overline Y&\stackrel{f_1}{\rightarrow}&\overline X\\
\downarrow&&\downarrow\\
Y&\stackrel{f}{\rightarrow}&X
\end{array}
$$
of the form of (\ref{eqI1}) such that $\overline Y$ is nonsingular, $f_1$ is ``close'' to
being finite, and the vertical arrows are birational and ``close'' to being proper.

Let $K$ be an algebraic function field over a field $k$ of characteristic zero. In his work on 
resolution of singularities, Zariski \cite{Z2} constructed for each $k$-valuation ring
$V$ of $K$ a projective model $X_V$ of $K$ such that the center of $V$ is nonsingular on $X_V$.
By the quasi-compactness of the Zariski-Riemann manifold of valuations of $K$, there exists
a finite number of the $X_V$, $\{X_{V_1},\ldots,X_{V_n}\}$ such that every valuation ring
$V$ of $K$ has a nonsingular center on at least one of the $X_{V_i}$.

In dimension $\le 3$, Zariski \cite{Z5} was able to patch the open nonsingular locus of 
appropriate birational transforms of the
$X_{V_i}$ to produce a nonsingular projective model $X$ of $K$.

The only part of Zariski's proof of the existence of a nonsingular model which does not
extend to arbitrary dimension is the final step where nonsingular open subsets $U_i$ of $X_{V_i}$ are patched (after appropriate birational transforms) to produce a projective
variety. Hironaka observes in Chapter 0, Section 6 of \cite{Ha} that we can always patch the
nonsingular loci $U_i$ of the $X_{V_i}$ along open sets where they are isomorphic, to produce
an integral finite type scheme $X$ such that $X$ is nonsingular and every valuation ring
$V$ of $K$ has a center on $X$, but $X$ will in general not be separated. Hironaka calls such
schemes ``complete''.

If such an $X$ is separated, then the following Lemma shows that $X$ is in fact complete in the  
usual sense, that is $X$ is a proper $k$-scheme.

\begin{Lemma}\label{LemmaS} Suppose that $X$ is a separated, integral, finite type $k$-scheme, such that
every $k$-valuation of $k(X)$ has a center on $X$. Then $X$ is a proper $k$-scheme.
\end{Lemma}

The proof of this Lemma is immediate from (ii) of Corollary II 7.3.10 \cite{EGA} (Recall that
in the notation of \cite{EGA},
a scheme is a separated pre-scheme), or can be deduced  directly from  the valuative criterion of properness (Theorem II 7.3.8 \cite{EGA} or Theorem II 4.7 \cite{Ha}).
 
One may expect a scheme $X$ which satisfies all of the conditions of the Lemma above except the
separatedness condition to be universally closed. This is false, as is shown by the following
example.

\begin{Example} Suppose that $k$ is a field. There exists a (nonseparated) integral finite type $k$-scheme such that
every $k$-valuation ring  of $k(X)$ has a center on $X$, but $X$ is not universally closed
over $k$.
\end{Example}

\begin{pf} Let $\phi$ be an imbedding of $\bold{P}^1$ in $\bold P^3$. Let $Z=\phi({\bold P}^1-\{\infty\})$, $x_0=\phi(\infty)$. Let $\pi:X_1\rightarrow{\bold P}^3$ be the blowup of
$\overline Z=\phi({\bold P}^1)$, $X_2={\bold P}^3-\{x_0\}$. We can construct an integral
finite type $k$-scheme $X$ by glueing $X_1$ to $X_2$ along the open sets $X_1-\pi^{-1}(\overline Z)$ and $X_2-Z$. By construction, every $k$-valuation of $k(X)$ has a center
on $X$. $\phi$ induces an isomorphism of ${\bold P}^1-\{\infty\}$ with the closed subset $Z\subset X_1$.
$Z$ is closed in $X$. Thus the induced morphism $\phi:\text{spec}(k({\bold P}^1)
\rightarrow X$ does not extend to a morphism $\text{spec}({\cal O}_{{\bold P}^1,\infty})
\rightarrow X$.

Suppose that $X$ is universally closed over $k$. Let $U=\text{spec}(k({\bold P}^1))$,
$T=\text{spec}({\cal O}_{{\bold P}^1,\infty})$, with natural morphism $i:U\rightarrow T$.
Let $t_1\in T$ be the generic point, $t_0\in T$ the special point.

Let $A$ be the closure of $\phi\times i(U)$ in $X\times T$. $\pi_2(A)$ is closed
since the projection $\pi_2:X\times T\rightarrow T$ is closed by assumption. So there exists $y_0\in A$
such that $\pi_2(y_0)=t_0$. By Lemma II 4.4 \cite{Ha} we have an extension
$T\rightarrow X\times T$ of $U\rightarrow X\times T$ which projects to an extension $\text{spec}(T)\rightarrow X$ of $U\rightarrow X$, and thus an extension of
$\phi$ to $\text{spec}({\cal O}_{{\bold P}^1,\infty})\rightarrow X$, a contradiction.
\end{pf}
\

\begin{Theorem}\label{Theorem8} Suppose that $f:Y\rightarrow X$ is a dominant, generically finite morphism
of complete varieties over a field $k$ of characteristic zero.  Then there exists a commutative diagram of integral, finite type
 $k$-schemes
$$
\begin{array}{lll}
\overline Y &\rightarrow &\overline X\\
\downarrow&&\downarrow\\
Y&\rightarrow&X
\end{array}
$$
such that $\overline Y$ is nonsingular, $\overline X$ has normal toric singularities, the vertical arrows are birational morphisms,
$\overline Y\rightarrow \overline X$ is quasi-finite and every $k$-valuation ring of
$k(X)$ has a center on $\overline X$, every $k$-valuation ring of $k(Y)$ has a center on $\overline Y$.
\end{Theorem}

\begin{Remark} Theorem \ref{Example1} (and Lemma \ref{LemmaS}) show that we cannot take the vertical arrows in Theorem \ref{Theorem8} to be proper. 
\end{Remark}

\begin{pf} Let $K=k(X)$ be the function field of $X$, $K^*=k(Y)$ be the function field
of $Y$. By assumption, $k(Y)$ is finite
over $k(X)$. By resolution of singularities \cite{H1}, we may assume 
that $Y$ and $X$ are nonsingular. Let $V^*$ be a $k$-valuation ring of $K^*$
such that $\text{trdeg}_kV^*/m_{V^*}=0$,
$V=V^*\cap K$. Let $p$ be the center of $V$ on $X$, $q$ the center of $V^*$ on $Y$.

By Theorem \ref{Theorem1}, there exists a sequence of the form (\ref{eq3}),
$$
\begin{array}{lllll}
R_0&\rightarrow R\rightarrow&S&\subset&V^*\\
\uparrow&&\uparrow\\
{\cal O}_{X,p}&\rightarrow&{\cal O}_{Y,q}
\end{array},
$$
such that $S$ is regular and $R$ has toric singularities.

Let $N={\bold Z}^n$, $\sigma$ be the cone generated by the rows of $A=(a_{ij})$ (with the
notation of (\ref{eq1*}) in $N\otimes{\bold R}$. Let $M$ be the dual lattice of $N$, $\hat\sigma$ be the dual cone of $\sigma$. By the proof of Theorem
\ref{Theorem1} (and Remark \ref{RemarkT}),  if
$\hat\sigma\cap M$ is generated by $e_1,\ldots, e_r$, then
$R=R_0[x^{e_1},\ldots,x^{e_r}]_m$.  There is a natural inclusion
$$
k[\hat\sigma\cap M]=
k[x^{e_1},\ldots,x^{e_r}]\rightarrow R.
$$
 $U_{\sigma}=\text{spec}(k[\hat\sigma\cap M])$ is a
normal affine toric
variety.

Thus there exist affine open sets $U_p$ of $p$ in $X$ and $\overline U_q$ of $q$ in $Y$ and
affine rings $R_V$ with quotient fields $K$ and $S_{V^*}$ with quotient field $K^*$ with
the following properties:
\begin{enumerate}
\item $R_V$ is normal and $S_{V^*}$ is regular.
\item If $p_1$ is the center of $V$ on $R_V$ and $q_1$ is the center of $V^*$ on $S_{V^*}$,
then $(R_V)_{p_1}=R$, $(S_{V^*})_{q_1}=S$.
\item There is a commutative diagram 
\begin{equation}\label{eq27}
\begin{array}{lll}
R_V&\rightarrow&S_{V^*}\\
\uparrow&&\uparrow\\
\Gamma(U_p,{\cal O}_X)&\rightarrow&\Gamma(\overline U_q,{\cal O}_Y)
\end{array}
\end{equation}
such that $R_V\rightarrow S_{V^*}$ is quasi-finite, 
$$
k[\hat\sigma\cap M]=k[x^{e_1},\ldots, x^{e_r}]\rightarrow R_V
$$
is \'etale, so that $R_V$ has normal toric singularites,
and the vertical arrows are birational
morphisms.
\end{enumerate}

Let $Z(K)$ be the Zariski-Riemann manifold of $k$-valuation rings of $K$, $Z(K^*)$
be the Zariski-Riemann manifold of $k$-valuation rings  of $K^*$. 
These spaces have natural topologies with respect to which they are quasi-compact (Theorem 40 Section 17, Chapter VI, \cite{ZS}).
     There is a natural continuous map
$\Phi:Z(K^*)\rightarrow Z(K)$ defined by $\Phi(V^*)=V^*\cap K$. 
For each $V^*\in Z(K^*)$ such that $\text{trdeg}_k V^*/m_{V^*}=0$, let $Y_{V^*}$ be a  projective variety which contains $\text{spec}(S_{V^*} )$ as an open set, and let $X_V$ be a  projective variety which
contains $\text{spec}(R_V)$ as an open set and such that the birational rational maps 
$Y_{V^*}\rightarrow X_V$, $Y_{V^*}\rightarrow Y$ and $X_V\rightarrow X$ are morphisms. Then there are commutative diagrams of
continuous maps
$$
\begin{array}{lll}
Z(K^*)&\rightarrow&Z(K)\\
\downarrow\pi_{V^*}&&\downarrow\pi_V\\
Y_{V^*}&\rightarrow&X_V\\
\downarrow&&\downarrow\\
Y&\rightarrow&X
\end{array}.
$$
Let $\overline W_{V^*}=\pi_{V^*}^{-1}(\text{spec}(S_{V^*}))$,  an open neighborhood of
$V^*$ in $Z(K^*)$,
$W_V=\pi_V^{-1}(\text{spec}(R_V))$, an open neighborhood of $V$ in $Z(K)$.
Suppose that $V'\in Z(K^*)$ and $\text{trdeg}_k V'/m_{V'}>0$. Let $\pi:V'\rightarrow
V'/m_{V'}$ be the residue map. By corollary 3 to Theorem 5, Section 4, Chapter VI \cite{ZS},
there exists a $k$-valuation ring $V_1$ with quotient field $V'/m_{V'}$ such that
$\text{trdeg}_k V_1/m_{V_1}=0$. Set $V^*=\pi^{-1}(V_1)$, a $k$-valuation ring of $K^*$
such that $V^*$ is a specialization of $V'$ and $\text{trdeg}_kV^*/m_{V^*}=0$ (page 57 \cite{Ab4}). Since $V^*$ dominates a local ring $S$ of $W_{V^*}$, $V'$ must dominate a localization of $S$, which is the local ring of a point of $W_{V^*}$. Thus $V'\in \overline W_{V^*}$.
$\{\overline W_{V^*}\}_{V^*\in Z(K^*)}$ is thus an open cover of $Z(K^*)$.
$\{W_{V}\}_{V\in Z(K)}$ is also an open cover of $Z(K)$.

Since $Z(K^*)$ and $Z(K)$ are quasi-compact, there is a finite set of valuation rings
$V_1^*,\ldots,V_m^*\in Z(K^*)$ with $\text{trdeg}_k V_i^*/m_{V_i^*}=0$ for all $i$ such that if $V_i=V_i^*\cap K$ then
$\{\overline W_{V_1^*},\ldots,\overline W_{V_m^*}\}$ is an open cover of $Z(K^*)$ and
$\{W_{V_1},\ldots,W_{V_m}\}$ is an open cover of $Z(K)$.

For $1\le i\le m$ we have commutative diagrams
$$
\begin{array}{rrr}
D_i:=\text{spec}(S_{V_i^*})&\rightarrow&C_i:=\text{spec}(R_{V_i})\\
b_i\downarrow&&a_i\downarrow\\
\overline U_{q_i}&\rightarrow&U_{p_i}
\end{array},
$$
where $q_i$ is the center of $V_i^*$ on $Y$, $p_i$ is the center of $V_i$ on $X$.
Let  $A_i\subset X$ be
nontrivial 
open sets where $a_i$ is an isomorphism and let $B_i\subset f^{-1}(A_i)$ be  nontrivial open sets where $b_i$ is an isomorphism. Then define $\overline Y$ to be the finite type
$k$-scheme obtained by patching $D_i$ to $D_j$ for $i\ne j$ along the nontrivial open
set $B_i\cap B_j$, and define $\overline X$ to be the finite type $k$-scheme obtained by patching $C_i$ to $C_j$ for $i\ne j$ along the nontrivial open set $A_i\cap A_j$.

By construction, $\overline Y$ and $\overline X$ are integral, $\overline Y$ is nonsingular, $\overline X$ has normal toric singularities and $\overline f:\overline Y\rightarrow \overline X$ is
quasi-finite.
\end{pf}

\section{Galois extensions}

For Galois extensions, Question 1 of the introduction has a positive answer. 
Suppose that $K\rightarrow K^*$ is a finite Galois extension of algebraic function fields
of characteristic zero. The existence of a finite map
of normal projective $k$-varieties $Y\rightarrow X$, where $Y$ is nonsingular,
$k(X)=K$, $k(Y)=K^*$  has been proven by Abhyankar in
Theorem 7 \cite{Ab3}.  We prove a relative version of this result in Theorem \ref{Theorem7}.

\begin{Theorem}\label{Theorem7} Suppose that $\Phi:Y\rightarrow X$ is a dominant morphism of projective
varieties over an algebraically closed field $k$ of characteristic zero such that $k(Y)$ is a finite
Galois extension of $k(X)$. Then there exists a commutative diagram
$$
\begin{array}{lll}
\overline Y&\rightarrow &\overline X\\
\downarrow&&\downarrow\\
Y&\rightarrow&X
\end{array}
$$
such that $\overline Y\rightarrow Y$,  $\overline X\rightarrow X$ are birational morphisms of projective
$k$-varieties, $\overline Y$ is
nonsingular, $\overline X$ is normal, $\overline Y\rightarrow \overline X$ is finite and $\overline X$ has normal toric singularities. 
\end{Theorem}

\begin{pf} Let $X_0$ be the normalization of $X$ in $k(X)$, $Y_0$ be the normalization of $X_0$
in $k(Y)$. Let $G=\text{Gal}(k(Y)/k(X))$. $G$ acts on $Y_0$ and $Y_0/G=X_0$. Let $D_0$
be the branch locus of $Y_0\rightarrow X_0$. Let $\pi_1:X_1\rightarrow X_0$ be a
resolution of singularities so that $D_1=\pi_1^{-1}(D_0)_{red}$ is a simple normal
crossings divisor. 
If $X_1$ is the blowup of an ideal sheaf ${\cal I}_0$ in $X_0$, let $Y_1$ be the normalization of the blowup of ${\cal I}_0{\cal O}_{Y_0}$.
$f_1:Y_1\rightarrow X_1$ is finite, and $G$ acts on $Y_1$. The branch locus 
of $f_1$ which is a divisor supported on $D_1$ (by the purity of the branch locus) has simple normal crossings, so by Abhyankar's Lemma (\cite{Ab7}, XIII 5.3 \cite{SGA}) $Y_1$ has 
normal toric singularities, and (by Lemma 7 \cite{Ab1}) if $p\in Y_1$ is a closed point, the stabilizer 
$$
G^s(p)=\{\sigma\in G\mid \sigma(p)=p\}
$$ 
is Abelian.
Suppose that ${\cal I}_1\subset {\cal O}_{Y_1}$ is an ideal sheaf such that the blowup of ${\cal I}_1$ in $Y_1$ dominates $Y$. Let ${\cal J}=\prod_{g\in G}g({\cal I}_1)$, and $\tilde
Y_2$ be the blowup of ${\cal J}$. Then $G$ acts on $\tilde Y_2$ and the rational map
$\tilde Y_2\rightarrow Y$ is a morphism. Let $Y_2\rightarrow \tilde Y_2$ be a
$G$-equivariant resolution of singularities of $\tilde Y_2$ (\cite{BM}, \cite{Vi}), with
composed map $\pi_2:Y_2\rightarrow Y_1$.
 $G$ acts on $Y_2$ and for $q\in Y_2$, $G^s(q)<G^s(\pi_2(q))$ so that $G^s(q)$ is Abelian. Let $\overline Y=Y_2$, $\overline X=Y_2/G$, a projective $k$-variety (c.f. page 126, \cite{Har}),
which has normal toric singularities (by Lemma 7, \cite{Ab1}).
\end{pf}

\begin{Example}\label{Example2}
Even if $k(Y)$ is Galois over $k(X)$, with Galois group $G$, and $Y\rightarrow X$ is
$G$-equivariant, we cannot take both $\overline X$ and $\overline Y$ to be nonsingular in Theorem \ref{Theorem8}
or in Theorem \ref{Theorem7}. 
\end{Example}

This is an immediate consequence of Abhyankar's example,
Theorem 11, \cite{Ab3} (restated in Theorem \ref{Theorem10} of this paper.)
Let $k$ be an algebraically closed field of characteristic zero. With the notations of
Theorem \ref{Theorem10}, $L_1$ is a Galois extension of $\overline K=k(u,v)$ with Galois group
${\bold Z}_q$. ${\bold Z}_q$ acts on $S=k[u,v,z]/z^q-uv^2$ and its invariant ring is
$R=k[u,v]$. With the notation of Theorem \ref{Theorem10}, $\overline R=R_{(u,v)}$.
By equivariant resolution of singularites \cite{BM}, \cite{Vi} (applied to the
normalization of $X={\bold P}^2$ in $L_1$) there exists a dominant 
${\bold Z}_q$ equivariant morphism of nonsingular
projective $k$-surfaces $Y\rightarrow X$ such that $k(X)=\overline K$, $k(Y)=L_1$ and there exists a point $p\in X$ such that ${\cal O}_{X,p}=\overline R$. 

Suppose that there exists a diagram 
$$
\begin{array}{lll}
\overline Y&\stackrel{\overline \Phi}{\rightarrow} &\overline X\\
\downarrow&&\downarrow\\
Y&\rightarrow&X
\end{array}
$$
as in the conclusions of Theorem \ref{Theorem8} such that $\overline X$ is nonsingular.
Let $\overline\nu$ be the valuation of $\overline K$ constructed in Theorem \ref{Theorem10}.
Let $\nu_1$ be the (unique) extension of $\overline\nu$ to $L_1$. Let $q\in\overline Y$ be a
center of $\nu_1$, $p=\overline\Phi(q)$. $B={\cal O}_{\overline Y,q}$ dominates 
$A={\cal O}_{\overline X,p}$. Since $\overline\Phi$ is quasi-finite, $B$ lies over $A$ by
Zariski's Main Theorem (10.9 \cite{Ab5}). Since $A$ is regular and $\overline\nu$ 
dominates $A$, $B$ is not regular by 
 Theorem \ref{Theorem10}, a contradiction.

\ \\ \\
\noindent
Steven Dale Cutkosky, Department of Mathematics, University of
Missouri\\
Columbia, MO 65211, USA\\
cutkoskys@@missouri.edu

\end{document}